\documentstyle[leqno,11pt]{article}

\newenvironment{dfntn}{\ \\
\noindent{\bf Definition.}\ \ }{\\ \ \par}

\newtheorem{thm}{Theorem}
\newtheorem{lemma}{Lemma}
\newtheorem{prop}{Proposition}

\newcommand{\contr}{\,\rule{.1in}{.5pt}\rule{.5pt}{1.5mm}\,\,}
\newcommand{\surj}{\,\longrightarrow\hspace{-12pt}
                     \longrightarrow\,\,}
\newcommand{\downsurj}{\ \makebox[0in]
     {\raisebox{-2pt}{$\downarrow$}}\makebox[0in]{$\downarrow$}\ }

\newcommand{\udot}{^{\textstyle\cdot}}
\newcommand{\pf}{\ \\ \noindent {\bf Proof.\ \ }}
\newcommand{\qed}{\ $\Box$\\ \ \par}
\newcommand{\rk}{\ \\ \noindent {\bf Remark.\ \ }}
\newcommand{\II}{{\rm II}}
\newcommand{\boxtimes}{\makebox[0in][l]{$\times$}\raisebox{-1pt}{$\Box$}}
\newcommand{\scriptboxtimes}{\makebox[0in][l]{$\scriptstyle\times$}
\makebox{$\scriptstyle\Box$\,}}
\newcommand{\displayboxtimes}
{\makebox[0in][l]{$\displaystyle\times$}\raisebox{-1pt}{$\displaystyle\Box$}}

\newcommand{\varpounds}{\makebox[-1.4pt][l]
 {\raisebox{1.7pt}{-}}{\cal L}}
\newcommand{\boldpounds}{\makebox[-1.4pt][l]
 {\raisebox{1.7pt}{\bf-}}{\mbox{\boldmath$\cal L$}}}

\title{Higher-order differentials of the period map and higher
Kodaira-Spencer classes}
\author{Yakov Karpishpan
\thanks{Supported in part by the NSF grant
DMS-9102233}}
\date{{\em Revised April} 1994}

\begin{document}
\maketitle

In \cite{K} we introduced two variants of higher-order
differentials of the period map and showed how to compute them for
a variation of Hodge structure that comes from  a deformation of a
compact K\"ahler manifold. More recently there appeared several
works (\cite{BG}, \cite{EV}, \cite{R}) defining higher tangent
spaces to the moduli and the corresponding higher Kodaira-Spencer
classes of a deformation. The $n^{th}$ such class $\kappa_n$
captures all essential information about the deformation up to
$n^{th}$ order.

A well-known result of Griffiths states that the (first)
differential of the period map depends only on the (first)
Kodaira-Spencer class of the deformation. In this paper we show
that the  second differential of the Archimedean period map
associated to a deformation is determined by $\kappa_2$ taken
modulo the image of $\kappa_1$, whereas the second differential
of the usual period map, as well as the second fundamental form
of the VHS, depend only on $\kappa_1$ (Theorems 2, 5, and 6 in
Section~3).

Presumably,  similar statements are valid in higher-order cases
(see Section~4).

\section{Constructing linear maps out of connections}

We start by reviewing the definitions of higher-order
differentials of the period map from \cite{K}, using a slightly
different approach. Let $S$ be a polydisc in ${\bf C}^s$
centered at 0. Consider a free ${\cal O}_S$-module $\cal V$ with
a decreasing filtration by ${\cal O}_S$-submodules
$\ldots\subseteq{\cal F}^p\subseteq{\cal
F}^{p-1}\subseteq\ldots$ and an integrable connection
$$\nabla:{\cal V}\rightarrow{\cal V}\otimes\Omega_S^1$$
satisfying the Griffiths transversality condition
$\nabla({\cal F}^p)\subseteq{\cal F}^{p-1}\otimes\Omega_S^1$.

\begin{lemma}
  (a) $\nabla$ induces an \ ${\cal O}_S$-linear map
  \begin{eqnarray*}
    \Theta_S & \longrightarrow & {\cal H}om_{{\cal O}_S} ({\cal
       F}^p,{\cal F}^{p-1}/{\cal F}^p)\\
    \xi & \longmapsto & \nabla_{\xi} \bmod {\cal F}^p
  \end{eqnarray*}

  (b) Analogously, we also have an \ ${\cal O}_S$-linear
symmetric map
   \begin{eqnarray*}
   \Theta_S^{\otimes 2} & \longrightarrow & {\cal H}om_{{\cal
         O}_S}
    ({\cal F}^p,{\cal F}^{p-2}/{\cal F}^p+span\,\{\nabla_{\eta}
    ({\cal F}^p)\,|\ all\ \eta\in\Gamma(S,\Theta_S)\}) \\
    \zeta\otimes\xi & \mapsto &
    \nabla_{\zeta}\nabla_{\xi} \bmod {\cal F}^p +
      span\,\{\nabla_{\eta}({\cal F}^p)\}\ .
    \end{eqnarray*}
\end{lemma}

\pf Both (a) and (b) are proved by straightforward computations;
the fact that the map in (b) is symmetric follows from the
integrability of $\nabla$:
$$\nabla_{\zeta}\nabla_{\xi}-
\nabla_{\xi}\nabla_{\zeta}=\nabla_{[\zeta,\xi]}\ ,$$
and so
$$\nabla_{\zeta}\nabla_{\xi}\equiv\nabla_{\xi}\nabla_{\zeta}
\bmod span\,\{\nabla_{\eta}({\cal F}^p)\,|\ \eta\in\Theta_S\}\
.$$
\qed

We will apply the Lemma to two connections arising from a
deformation of a compact K\"{a}hler manifold $X$,
\begin{eqnarray}
    \nonumber {\cal X} & \supset & X\\
    \pi\ \downarrow &   & \downarrow\\
    \nonumber S & \ni & 0
\end{eqnarray}

\noindent 1. The usual Gauss-Manin connection $\nabla$ on
${\cal H}=R^m\pi_{\ast}{\bf C}_{\cal X}$\ . In this case we
denote the map in (a)
$$d\Phi: \Theta_S\rightarrow
{\cal H}om_{{\cal O}_S} ({\cal F}^p,{\cal F}^{p-1}/{\cal F}^p)\
. $$
This is {\em the (first) differential of the (usual) period
map}. The same notation and terminology will be applied to the
induced map  $$\Theta_S\rightarrow\bigoplus_p{\cal H}om_{{\cal
O}_S}
  ({\cal F}^p/{\cal F}^{p+1},{\cal F}^{p-1}/{\cal F}^p)\ .$$

The map given by part (b) of the Lemma is called {\em the
second differential of the (usual) period map} and denoted
$$d^2\Phi: \Theta_S^2\longrightarrow{\cal H}om_{{\cal O}_S}
({\cal F}^p,{\cal F}^{p-2}/{\cal F}^p+span\,\{\nabla_{\eta}
    ({\cal F}^p)\,|\ \eta\in\Theta_S\})\ .$$

\noindent 2. The Archimedean Gauss-Manin connection
$\nabla=\nabla_{ar}$ on
$${\cal H}\otimes B_{ar}= R^m\pi_{\ast}{\bf
C}_{\cal X}[T,T^{-1}]$$
(see Appendix). The corresponding map from (a) is called
{\em the (first) differential of the {\em Archimedean}
period map}, denoted $$d\Psi:\Theta_S\longrightarrow
{\cal H}om_{{\cal O}_S}({\cal H}_{ar},
{\cal F}_{ar}^{-1}/{\cal H}_{ar})\ .$$
Again, we will abuse notation and write $d\Psi$ to denote the
induced map
\begin{equation}
\Theta_S\rightarrow {\cal H}om_{{\cal O}_S}
(Gr_{{\cal F}_{ar}}^0, Gr_{{\cal F}_{ar}}^{-1})\ .
\end{equation}
Finally, the map in part (b) of the Lemma, for
$\nabla=\nabla_{ar}$, is {\em the second differential of the
Archimedean period map} and will be denoted
$$d^2\Psi:\Theta_S^2\longrightarrow {\cal H}om_{{\cal O}_S}
({\cal H}_{ar},{\cal F}_{ar}^{-1}/{\cal H}_{ar}+
span\,\{\nabla_{\eta}({\cal H}_{ar})\,|\ \eta\in\Theta_S\})\ .$$

We have an identification of ${\cal O}_S$-modules
$$Gr_{{\cal F}_{ar}}^0({\cal H}\otimes B_{ar})\cong
{\bf R}^m\pi_{\ast}(Gr_F^0(\Omega\udot_{{\cal X}/S}\otimes
B_{ar}))\cong{\bf R}^m\pi_{\ast}\Omega\udot_{{\cal X}/S}=
{\cal H}\ ,$$
obtained from the obvious isomorphism of sheaf complexes

$$\begin{array}{ccccc}
  \longrightarrow & \Omega^p_{{\cal X}/S}.T^p &
  \stackrel{\bf d}{\longrightarrow} &
  \Omega^{p+1}_{{\cal X}/S}.T^{p+1} & \longrightarrow\\
    & \simeq\downarrow & & \downarrow\simeq & \\
  \longrightarrow & \Omega^p_{{\cal X}/S} &
  \stackrel{d}{\longrightarrow} &
  \Omega^{p+1}_{{\cal X}/S} & \longrightarrow
\end{array}$$

\noindent (``dropping the T's").

Similarly, $Gr_{{\cal F}_{ar}}^{-1}\cong {\cal H}$.
We use these identifications to obtain a version of (2)
``without
T's," $$\overline{d\Psi}:\Theta_S\longrightarrow{\cal E}nd_
{{\cal O}_S}({\cal H})\ .$$

Analogously, we also define
$$\overline{d^2\Psi}:\Theta_S^{\otimes 2}\longrightarrow{\cal
E}nd_ {{\cal O}_S}({\cal H})\ ,$$
with $\overline{d^2\Psi}(\zeta,\xi)$ being the composition
$${\cal H}\cong Gr_{{\cal F}_{ar}}^0
\stackrel{d^2\Psi(\zeta,\xi)}{\longrightarrow}{\cal
F}_{ar}^{-2}/{\cal H}_{ar}+span\,\{\nabla_{\eta}({\cal
F}_{ar}^{-1})\,|\ \eta\in\Theta_S\}\surj
Gr_{{\cal F}_{ar}}^{-2}\cong{\cal H}\ .$$

\rk Let ${\bf t}=(t_1,\ldots,t_s)$ be a coordinate system on
$S$ centered at 0. Then, in the
notation of \cite{K}, $d\Psi(\partial/\partial t_i)|_{t=0}$
is $\widetilde{L}^{(i)}$,
$\overline{d\Psi}(\partial/\partial t_i)|_{t=0}$ is
$\overline{L}^{(i)}$, $d^2\Psi(\partial/\partial t_i\otimes
\partial/\partial t_j)|_{t=0}$ is $\widetilde{L}^{(ij)}$, and
$\overline{d^2\Psi}(\partial/\partial t_i\otimes
\partial/\partial t_j)|_{t=0}$ is $\overline{L}^{(ij)}$.
\ ${\cal O}_S$-linearity of $d\Psi$, $d^2\Psi$, etc. is essential
for the ability to restrict to 0 in $S$.

To formulate the next Lemma, we bring out the natural
$C_S^{\infty}$-linear identification
$$h:{\cal H}=\bigoplus_{p+q=m}{\cal H}^{p,q}
\stackrel{\simeq}{\longrightarrow}
\bigoplus_p Gr_{\cal F}^p{\cal H}=:
Gr_{\cal F}\udot{\cal H}\ .$$

\begin{lemma}
  (a) For any \ $\xi\in\Theta_S$ and all \ $p$ we have
  $$\overline{d\Psi}(\xi)({\cal F}^p{\cal H})\subset
  {\cal F}^{p-1}{\cal H}\ ,$$
  the induced endomorphism of degree  \
  $-1$ of \ $Gr_{\cal F}\udot{\cal H}$ coincides with \ $d\Phi(\xi)$
  and, in fact,
  $$\overline{d\Psi}(\xi)=h^{-1}\circ d\Phi(\xi)\circ h\ .$$

  (b) For any \ $\zeta,\xi\in\Theta_S$ and all \ $p$ we have
  $$\overline{d^2\Psi}(\zeta,\xi)({\cal F}^p{\cal H})\subset
  {\cal F}^{p-2}{\cal H}\ ,$$ and
  $$d^2\Phi(\zeta,\xi):{\cal F}^p\longrightarrow{\cal F}^{p-2}/
  {\cal F}^p+span\,\{\nabla_{\eta}({\cal F}^p)\,|\
  \eta\in\Theta_S\}$$
  factors through \ $\overline{d^2\Psi}(\zeta,\xi):
  {\cal F}^p\longrightarrow{\cal F}^{p-2}$.
\end{lemma}

\pf (a) For any $\xi\in\Theta_S$
$$\nabla_{\xi}{\cal H}^{p,q}\subset
{\cal H}^{p,q}\oplus{\cal H}^{p-1,q+1}$$
and, correspondingly,
$$\nabla_{\xi}^{ar}{\cal H}^{p,q}.T^p\subset
{\cal H}^{p,q}.T^p\oplus{\cal H}^{p-1,q+1}.T^p\ .$$
Therefore, $d\Psi(\xi)$ maps ${\cal H}^{p,q}.T^p$ to its image
under $\nabla_{\xi}^{ar}$ modulo ${\cal H}_{ar}={\cal
F}_{ar}^0$, i.e. into ${\cal H}^{p-1,q+1}.T^p$. Hence
$$\overline{d\Psi}(\xi)({\cal H}^{p,q})\subset
{\cal H}^{p-1,q+1}\ ,$$
which implies every statement in part (a) of the Lemma.

Part (b) is established by similar reasoning.\qed

The connection $\nabla$ on $\cal H$ naturally induces a connection
$\nabla$ on ${\cal E}nd_{{\cal O}_S}({\cal H})$ subject to the
rule $$\nabla(Ax)=(\nabla A)x+A\nabla x$$
for any $A\in{\cal E}nd_{{\cal O}_S}({\cal H})$ and $x\in{\cal H}$.

In accordance with Lemma 1, for any $\zeta\in\Theta_S$
the covariant derivative $\nabla_{\zeta}$ on ${\cal E}nd({\cal
H})$ determines  an ${\cal O}_S$-linear map
$${\cal E}_{\zeta}:
im\,(\overline{d\Psi})\longrightarrow
		{\cal E}nd({\cal H})/im\,(\overline{d\Psi})\ .$$

\begin{lemma}
		For any \ $\zeta,\xi\in\Theta_S$

  (a) \ $\nabla_{\zeta}(\overline{d\Psi}(\xi))({\cal F}^p)\subset
  					{\cal F}^{p-1}$ for all \ $p$.

  (b) \ ${\cal E}_{\zeta}(\overline{d\Psi}(\xi))$ determines an
       element of \ ${\cal E}nd_{{\cal O}_S}^{(-1)}(Gr\udot
       _{\cal F}{\cal H})/im\,(d\Phi)$.

\end{lemma}

\pf (a) $\nabla_{\zeta}(\overline{d\Psi}(\xi))\omega=
           \nabla_{\zeta}(\overline{d\Psi}(\xi)\omega)-
           \overline{d\Psi}(\xi)\nabla_{\zeta}\omega$ for any
$\omega\in{\cal H}$. Assume $\omega\in{\cal H}^{p,q}$. We want to
show that the $(p-2,q+2)$-component of the right-hand side is 0.
But this component is
$$(\nabla_{\zeta}\nabla_{\xi}\omega)_{(p-2,q+2)}-
   (\nabla_{\xi}\nabla_{\zeta}\omega)_{(p-2,q+2)}=
   (\nabla_{[\xi,\zeta]}\omega)_{(p-2,q+2)}=0\ !$$

(b) follows from (a) and the relation between
$\overline{d\Psi}(\xi)$ and $d\Phi(\xi)$ established in part (a)
of the previous Lemma.\qed

\begin{dfntn}
		The {\em second fundamental form of the VHS}
  $$\II:\Theta_S^{\otimes 2}\longrightarrow
     {\cal E}nd_{{\cal O}_S}(Gr\udot_{\cal F})/im\,(d\Phi)$$
  is defined by
  $$\II(\zeta,\xi):=h\circ\nabla_{\zeta}(h^{-1}d\Phi(\xi)h)\circ
     h^{-1}\bmod im\,(d\Phi)\ .$$
\end{dfntn}

\rk $\II|_{t=0}$ was denoted $d^2\Phi$ in \cite{K}.

\

We will omit the identification $h$ in what follows.

\begin{prop}[(2.4) in \cite{K}]
$$\II(\zeta,\xi)\equiv\overline{d^2\Psi}(\zeta,\xi)-
   \overline{d\Psi}(\zeta) \circ \overline{d\Psi}(\xi)
   \bmod im\,(d\Phi)\ .$$

\end{prop}

\pf
$\nabla_{\zeta}(\overline{d\Psi}(\xi))\omega=
\nabla_{\zeta}(\overline{d\Psi}(\xi)\omega)-
\overline{d\Psi}(\xi)\nabla_{\zeta}\omega$ for any
$\omega\in{\cal H}$. Now, let $\mbox{\boldmath $\omega$}$ be
the element of $Gr_{{\cal F}_{ar}}^0({\cal H}\otimes B_{ar})$
corresponding to $\omega$ under the
isomorphism
${\cal H}\cong Gr_{{\cal F}_{ar}}^0({\cal H}\otimes B_{ar})$.
Using a similar identification of ${\cal H}$ with
$Gr_{{\cal F}_{ar}}^{-2}({\cal H}\otimes B_{ar})$, we have the
following correspondences:
$$\nabla_{\zeta}(\overline{d\Psi}(\xi)\omega)\longleftrightarrow
(\nabla_{\zeta}(d\Psi(\xi)\mbox{\boldmath$\omega$})
  \bmod {\cal F}_{ar}^{-1})=(d^2\Psi(\zeta,\xi)
  \mbox{\boldmath$\omega$}\bmod {\cal F}_{ar}^{-1})
  \in Gr_{{\cal F}_{ar}}^{-2}\ ,$$
and
$$\overline{d\Psi}(\xi)\nabla_{\zeta}\omega\longleftrightarrow
(d\Psi(\zeta)\circ d\Psi(\xi)\mbox{\boldmath$\omega$}
\bmod {\cal F}_{ar}^{-1})\in Gr_{{\cal F}_{ar}}^{-2}\ .$$
It remains to pass to ${\cal H}$ on the right-hand side, i.e.
put bars on $d\Psi$ and $d^2\Psi$.\qed

\section{The second Kodaira-Spencer class}

First, let us recall the construction of the (first)
Kodaira-Spencer map $\kappa_1=\kappa$ of the deformation (1): it
is the connecting morphism in the higher-direct-image sequence
\begin{equation}
  \longrightarrow\pi_{*}\Theta_{\cal
       X}\longrightarrow\Theta_S
       \stackrel{\kappa}{\longrightarrow}R^1\pi_{*}\Theta_{{\cal
       X}/S} \longrightarrow
\end{equation}
associated with the short exact sequence
\begin{equation}
   0\longrightarrow\Theta_{{\cal X}/S}\longrightarrow
    \Theta_{\cal X}\longrightarrow\pi^{*}\Theta_S\longrightarrow
        0\     .
\end{equation}

Given $\xi\in\Theta_S$, the corresponding covariant derivative of
the Gauss-Manin connection
$$\nabla_{\xi}:{\bf R}^m\pi_{*}\Omega\udot_{{\cal X}/S}
\longrightarrow{\bf R}^m\pi_{*}\Omega\udot_{{\cal X}/S}$$
is computed as follows (see \cite{Del}, \cite{KO}, or
\cite{K}). Choose a Stein covering ${\cal U}=\{U_i\}$ of $X$.
Then $\{W_i=U_i\times S\}$ constitute a Stein covering $\cal W$
of $\cal X$.  Consider a class in  $\Gamma(S,{\bf
R}^m\pi_{*}\Omega\udot_{{\cal X}/S})$ represented by the \v{C}ech
cocycle $\omega=\{\omega_Q\in \Gamma({\cal W}_Q,\Omega^p_{{\cal
X}/S})\}$, where $Q=(i_1<\ldots<i_q)$ and $p+q+1=m$. Let the same
symbol  $\omega_Q$ denote a pull-back of $\omega_Q\in
\Gamma({\cal W}_Q,\Omega^p_{{\cal X}/S})$ to an element of
$\Gamma({\cal W}_Q,\Omega^p_{\cal X})$.
Let $v=\{v_i\}$ denote liftings of $\xi\in\Theta_S$, or rather,
in $\pi^{*}\Theta_S$, to $\Gamma(W_i,\Theta_{\cal X})$. Then
$$\nabla_{\xi}[\omega]=[\check{\pounds}_v\omega]$$
where the brackets denote cohomology classes in
${\bf R}^m\pi_{*}\Omega\udot_{{\cal X}/S}$ (more precisely, in
$\Gamma(S,{\bf R}^m\pi_{*}\Omega\udot_{{\cal X}/S})$), and
$\check{\pounds}_v$ is the Lie derivative on
$\check{C}^m({\cal W},\Omega\udot_{\cal X}) $  with respect to
$v=\{v_i\}\in \check{C}^0({\cal W},\Theta_{\cal X})$:
\begin{equation}
		\check{\pounds}_v\omega=
     \{\check{D}v_{i_1}\contr\omega_{i_1,\ldots,i_q}+
     v_{i_1}\contr\check{D}\omega_{i_1,\ldots,i_q}\}
\end{equation}
with $\check{D}=d\pm\delta$, $\delta=\check{\delta}$ being the
\v{C}ech differential, as usual.

Now, when $\xi\in\Theta_S$ lifts to all of $\cal X$, i.e. $\xi$
lies in the image of $\pi_*\Theta_{\cal X}\rightarrow\Theta_S$
($= ker\,(\kappa)$ !), the cochain
$v\in\check{C}^0({\cal W},\Theta_{\cal X})$ lifting $\xi$ is a
cocycle, i,e, $\delta v=0$. But then  formula (5) reduces to
$$\check{\pounds}_v\omega=
\{dv_{i_1}\contr\omega_{i_1,\ldots,i_q}+
 v_{i_1}\contr d\omega_{i_1,\ldots,i_q}\}=
 \{\pounds_{v_{i_1}}\omega_{i_1,\ldots,i_q}\} ,$$
where $\pounds_{v_{i_1}}$ now denotes the usual Lie derivative
with respect to the vector field $v_{i_1}$. Evidently, in this
case $\nabla_{\xi}{\cal F}^p\subset{\cal F}^p$. Consequently,
the (first) differential of the period map
$$d\Phi:\Theta_S\longrightarrow
{\cal H}om({\cal F}^p,{\cal H}/{\cal F}^p)$$
factors through
$\Theta_S/\,im\,\{\pi_*\Theta_{\cal X}\rightarrow\Theta_S\}=
\Theta_S/\,ker\,(\kappa)$, and thus we arrive at

\begin{thm}[Griffiths]
There is a commutative
diagram
$$\begin{array}{ccccc}
			\Theta_S & \stackrel{d\Phi}{\longrightarrow} &
       \bigoplus_p{\cal H}om({\cal F}^p,{\cal H}/{\cal F}^p) &
       \surj & {\cal E}nd^{(-1)}(Gr\udot_{\cal F}{\cal H}) \\
     & & & & \\
     & \searrow\kappa & \uparrow & \smile\nearrow & \\
     & & & & \\
     &      & {\bf T}^1_{{\cal X}/S} &  &
\end{array}$$
where \ ${\bf T}^1_{{\cal X}/S}:=R^1\pi_*\Theta_{{\cal X}/S}$ and
the northeast arrow sends \ $x$ to the map \ $x\smile$ \ \ (the
cup product with \ $x$).
\end{thm}

Analogous results hold for $d\Psi$ and $\overline{d\Psi}$
(see \cite{K}). We seek a similar statement for $d^2\Phi$ and
$d^2\Psi$.

First we need to review the construction of the ``second-order
tangent space to the moduli." Here we are following (a
relativized version of) the presentation in \cite{R}. We will
make the simplifying assumption that $X$ has no global
holomorphic vector fields, i.e. $\pi_*\Theta_{{\cal X}/S}=0$.

Consider the diagram
\begin{eqnarray}
   \nonumber{\cal X}\times_S{\cal X} &
       \stackrel{g}{\longrightarrow} & {\cal X}_2/S \\
    f\searrow & & \swarrow p \\
    \nonumber & S &
\end{eqnarray}
where ${\cal X}_2/S$ denotes the symmetric
product of $\cal X$ with itself over $S$ (fiberwise).
Let ${\cal K}\udot$ denote the complex of sheaves on ${\cal X}_2/S$
$$\begin{array}{ccc}
  {\scriptstyle -1} & & {\scriptstyle 0} \\
  (g_*(\Theta_{{\cal X}/S}^{\scriptboxtimes 2}))^{-} &
     \stackrel{[\, ,\,]}{\longrightarrow} & \Theta_{{\cal X}/S}
\end{array}$$
where $\boxtimes$ stands for the exterior tensor product on
${\cal X}\times_S{\cal X}$, \ \ $(\ \ )^{-}$ denotes the
anti-invariants of the
${\bf Z}/2{\bf Z}$-action, and the differential is the
restriction to the diagonal
$\Delta\subset{\cal X}\times_{S}{\cal X}$,
followed by the Lie bracket of vector fields.

\begin{dfntn}
${\bf T}^{(2)}_{{\cal X}/S}:={\bf R}^1p_*{\cal K}\udot$ is
the sheaf on $S$ whose fiber over each $t\in S$ is the
{\em second-order (Zariski) tangent space to the base $V_t$ of
the miniversal deformation of} $X_t$, i.e.
${\bf T}^{(2)}_{X_t}\cong({\bf m}_{V_t,0}/{\bf m}_{V_t,0}^3)^*$.
\end{dfntn}

This should not be confused with the sheaf
${\bf T}^{2}_{{\cal X}/S}=R^2\pi_*\Theta_{{\cal X}/S}$ whose
fiber over each  $t\in S$ is the {\em obstruction space}\ \
${\bf T}^{2}_{X_t}$ for deformations of $X_t$.
In fact, we have this exact sequence:
\begin{equation}
		0\longrightarrow{\bf T}^1_{{\cal X}/S}\longrightarrow
    {\bf T}^{(2)}_{{\cal X}/S}\longrightarrow
    Sym^2{\bf T}^1_{{\cal X}/S}\stackrel{o}{\longrightarrow}
    {\bf T}^{2}_{{\cal X}/S}\ ,
\end{equation}
where $o$ is the first obstruction map, given by the
${\cal O}_S$-linear graded Lie bracket:
$$Sym^2R^1\pi_*\Theta_{{\cal
X}/S}\stackrel{[\, ,\,]}{\longrightarrow}
R^2\pi_*\Theta_{{\cal X}/S}\ .$$

We will find it easier to deal with an ``unsymmetrized version"
of ${\bf T}^{(2)}_{{\cal X}/S}$.

\begin{dfntn}
$\widetilde{\bf T}^{(2)}_{{\cal X}/S}:={\bf R}^1f_*
\widetilde{\cal K}\udot$,
where $\widetilde{\cal K}\udot$ is the complex on
${\cal X}\times_S{\cal X}$,
$$\begin{array}{ccc}
  {\scriptstyle -1} & & {\scriptstyle 0} \\
  \Theta_{{\cal X}/S}^{\scriptboxtimes 2} &
     \stackrel{[\, ,\,]}{\longrightarrow} & \Theta_{{\cal X}/S}\ .
\end{array}$$
\end{dfntn}
$\widetilde{\bf T}^{(2)}_{{\cal X}/S}$ fits in the commutative diagram
with exact rows:
$$\begin{array}{cccccccc}
  0\rightarrow & {\bf T}^1_{{\cal X}/S} & \longrightarrow &
    \widetilde{\bf T}^{(2)}_{{\cal X}/S} & \longrightarrow &
    ({\bf T}^1_{{\cal X}/S})^{\otimes 2} &
    \stackrel{\widetilde{o}}{\longrightarrow} &
    {\bf T}^2_{{\cal X}/S} \\
   & \| & & \downsurj & & \downsurj & & \| \\
  0\rightarrow & {\bf T}^1_{{\cal X}/S} & \longrightarrow &
    {\bf T}^{(2)}_{{\cal X}/S} & \longrightarrow &
    Sym^2{\bf T}^1_{{\cal X}/S} &
    \stackrel{o}{\longrightarrow} &
    {\bf T}^2_{{\cal X}/S}
\end{array}$$

\begin{dfntn}
 $T_S^{(2)}:={\cal D}_S^{(2)}/{\cal O}_S$ will denote the sheaf of
 {\em second-order tangent  vectors}\ \  on $S$.
\end{dfntn}

As part of a more general construction in \cite{EV}, there is
{\em the second Kodaira-Spencer map} associated to every
deformation as in (1):
$$\kappa_2:T_S^{(2)}\longrightarrow{\bf T}^{(2)}_{{\cal X}/S}\ .$$
We will work with a natural lifting $\widetilde{\kappa}_2$
of $\kappa_2$:
\begin{equation}
\begin{array}{ccc}
   \Theta_S\oplus\Theta_S^{\otimes 2} &
         \stackrel{\widetilde{\kappa}_2}{\longrightarrow} &
         \widetilde{\bf T}^{(2)}_{{\cal X}/S} \\
      \downsurj & & \downsurj \\
    T_S^{(2)} & \stackrel{\kappa_2}{\longrightarrow} &
          {\bf T}^{(2)}_{{\cal X}/S}
\end{array}
\end{equation}
or, rather, with the restriction of $\widetilde{\kappa}_2$ to
$\Theta_S^{\otimes 2}$.

It is easy to describe $\widetilde{\kappa}_2$ explicitly.
Let
$$\kappa:\Theta_S\longrightarrow {\bf T}^1_{{\cal X}/S}=
R^1\pi_*\Theta_{{\cal X}/S}$$
be the ({\em relative, first}\ ) Kodaira-Spencer map of the
family $\pi$ as in (1).
It is equivalent to the datum of a section of
$\Gamma(S,\Omega_S^1\otimes R^1\pi_*\Theta_{{\cal X}/S})$. This
section can be represented by a
$\check{C}^1({\cal U},\Theta_X)$-valued one-form on $S$,
\begin{equation}
  \theta({\bf t})d{\bf t}:=
  \sum_{\ell=1}^s\theta({\bf t})_{\ell}dt_{\ell}=
  \sum_{\ell=1}^s\sum_{I\in{\bf Z}_+^s}^s
    \theta_{\ell}^{(I)}{\bf t}^Idt_{\ell}\ \ \
    ({\bf Z}_+:=\{0\}\cup{\bf N})\ .
\end{equation}
Here each $\theta_{\ell}^{(I)}=\{\theta_{ij,\ell}^{(I)}\}_{ij}$ is
a cochain in $\check{C}^1({\cal U},\Theta_X)$ and each
$\theta({\bf t})_{\ell}$ is a cocycle on $X_t$ for every
value of $\bf t$, but only the leading coefficients
$\theta_{\ell}^{(0)}\ \ (\ell=1,\ldots,s)$ are \v{C}ech {\em
cocycles}\ on $X$\ \ $(t=0)$. The rest satisfy the ``deformation
equation"
\begin{equation}
  \delta\left(\frac{\partial\theta({\bf t})_{\ell}}
    {\partial_{t_k}}\right)=
  [\theta({\bf t})_{\ell},\theta({\bf t})_k]
\end{equation}
When $s=1$, this
reduces to $\delta\dot{\theta}(t)=[\theta(t),\theta(t)]$.

\noindent Thus, it is natural to make the following

\begin{dfntn}
 $\widetilde{\kappa}_2:\Theta_S^{\otimes 2}\longrightarrow
 \widetilde{\bf T}^{(2)}_{{\cal X}/S}$
 sends
 $\frac{\partial}{\partial t_k}\otimes
 \frac{\partial}{\partial t_{\ell}}$  to the cohomology class of
 the cocycle
\begin{equation}
 (\theta({\bf t})_k\times
 \theta({\bf t})_{\ell},
 \frac{\partial\theta({\bf t})_{\ell}}
    {\partial_{t_k}})
    \in \check{C}^1({\cal W}\times_S{\cal W},
    \widetilde{\cal K}\udot)\ .
\end{equation}
\end{dfntn}
For example, if
\begin{equation}
 \theta({\bf t})d{\bf t}=\sum_{\ell=1}^s(\theta_{\ell}^{(0)}+
 \sum_{k=1}^s\theta_{\ell}^{(k)}t_k)dt_{\ell}+O({\bf t}^2)
\end{equation}
is the expansion of $\theta({\bf t})d{\bf t}$ at 0 up to order
two, then
$\widetilde{\kappa}_2|_{0}:\Theta_S^{\otimes 2}|_{0}\rightarrow
\widetilde{\bf T}_X^{(2)}$ sends
$\frac{\partial}{\partial t_k}\otimes
\frac{\partial}{\partial t_{\ell}}$  to the cohomology class of
the cocycle
\begin{equation}
 (\theta_k^{(0)}\times\theta_{\ell}^{(0)},
 \theta_{\ell}^{(k)})\in \check{C}^1({\cal U}\times{\cal U},
 \widetilde{\cal K}\udot|_{0})\ .
\end{equation}

Indeed, for the definition of $\widetilde{\kappa}_2$ to make any
sense, we must have the following commutative diagram with exact
rows:
\begin{equation}
\begin{array}{cccccccc}
  0\rightarrow & \Theta_S & \rightarrow &
    \Theta_S\oplus\Theta_S^{\otimes 2} &
    \rightarrow & \Theta_S^{\otimes 2} & \rightarrow  & 0 \\
   & \kappa_1\ \downarrow & & \downarrow \widetilde{\kappa}_2 &
   & \downarrow \kappa_1^2 & & \\
  0\rightarrow & {\bf T}_{{\cal X}/S}^1 & \rightarrow &
    \widetilde{\bf T}_{{\cal X}/S}^{(2)} & \rightarrow &
    ({\bf T}_{{\cal X}/S}^1)^{\otimes 2} &
    \stackrel{o}{\rightarrow} & {\bf T}_{{\cal X}/S}^2\ .
\end{array}
\end{equation}
The square on the right induces a commutative triangle
\begin{equation}
\begin{array}{ccc}
  \Theta_S^{\otimes 2} &
   \stackrel{\widetilde{\kappa}_2}{\longrightarrow} &
   \widetilde{\bf T}_{{\cal X}/S}^{(2)} \\
   & \kappa_1^2\ \searrow & \downarrow \\
   & &  ({\bf T}_{{\cal X}/S}^1)^{\otimes 2}\ .
\end{array}
\end{equation}
Therefore,
$\widetilde{\kappa}_2(\frac{\partial}{\partial t_k}\otimes
\frac{\partial}{\partial t_\ell})$ must project onto
$$\kappa_1(\frac{\partial}{\partial t_k})\otimes
\kappa_1(\frac{\partial}{\partial t_\ell})=
[\theta({\bf t})_k]\otimes[\theta({\bf t})_\ell]\ .$$
Since
$$\check{C}^1({\cal W}\times_S{\cal W},\widetilde{\cal K}\udot)=
\check{C}^2({\cal W}\times_S{\cal W},\widetilde{\cal K}^{-1})
\oplus\check{C}^1(\widetilde{\cal K}^0)\ ,$$
and
$$\check{C}^2({\cal W}\times_S{\cal W},\widetilde{\cal K}^{-1})
\simeq\check{C}^1({\cal U},\Theta_{{\cal X}/S})^{\otimes 2}\ ,$$
this means that the $\check{C}^2(\widetilde{\cal
K}^{-1})$-component of a representative of
$\widetilde{\kappa}_2(\frac{\partial}{\partial t_k}\otimes
\frac{\partial}{\partial t_\ell})$ in
$\check{C}^1(\widetilde{\cal K}\udot)$ may be taken to be
$$\theta({\bf t})_k\times\theta({\bf t})_{\ell}\ .$$
And, in view of (10), the cochain (11) is indeed a {\em cocycle}
\  in  $\check{C}^1(\widetilde{\cal K}\udot)$.

We still need to check that $\widetilde{\kappa}_2$ is well-defined.

\begin{prop}\ \ \
$\widetilde{\kappa}_2:\Theta_S^{\otimes 2}\longrightarrow
\widetilde{\bf T}^{(2)}_{{\cal X}/S}$ can
be presented as a connecting morphism in the higher-direct-image
sequence of a short exact sequence.
\end{prop}

\pf  The starting point is the sequence of
${\cal O}_{\cal X}$-modules
(4), whose direct-image sequence (3) gives the first
Kodaira-Spencer mapping $\kappa_1$ as a connecting morphism.
Now, (4) contains an exact subsequence
\begin{equation}
 0\longrightarrow\Theta_{{\cal X}/S}
   \stackrel{\alpha}{\longrightarrow}\widetilde{\Theta}_{\cal X}
   \stackrel{\beta}{\longrightarrow}\pi^{-1}\Theta_S
   \longrightarrow 0
\end{equation}
of $\pi^{-1}{\cal O}_S$-modules, whose direct-image sequence also
has $\kappa_1$ as a connecting morphism. From now on we will
work with (16) in place of (4).

We can splice two sequences produced from (16)  by exterior
tensor products with $\Theta_{{\cal X}/S}$ and with
$\pi^{-1}\Theta_S$, respectively:
$$\begin{array}{cccc}
  & & 0 & \\
  & & \uparrow & \\
  & & \pi^{-1}\Theta_S\displayboxtimes\pi^{-1}\Theta_S & \\
  & & \uparrow & \\
  & & \pi^{-1}\Theta_S\displayboxtimes
    \widetilde{\Theta}_{\cal X} & \\
  & \nearrow & \uparrow & \\
 0\longrightarrow\Theta_{{\cal X}/S}\displayboxtimes
    \Theta_{{\cal X}/S} \longrightarrow
    \widetilde{\Theta}_{\cal X}\displayboxtimes
    \Theta_{{\cal X}/S} & \rightarrow &
    \pi^{-1}\Theta_S\displayboxtimes\Theta_{{\cal X}/S} &
    \longrightarrow 0 \\
  & & \uparrow & \\
  & & 0 &
\end{array}$$
The resulting four-term exact sequence
\begin{equation}
  0\rightarrow\Theta_{{\cal X}/S}^{\scriptboxtimes 2}
  \rightarrow\widetilde{\Theta}_{\cal X}
  \displayboxtimes\Theta_{{\cal X}/S}
  \stackrel{\beta\scriptboxtimes\alpha}{\longrightarrow}
  \pi^{-1}\Theta_S\displayboxtimes\widetilde{\Theta}_{\cal X}
  \rightarrow(\pi^{-1}\Theta_S)^{\scriptboxtimes 2}
  \rightarrow 0
\end{equation}
can be extended to a commutative diagram
\begin{equation}
\begin{array}{ccccc}
0\rightarrow & \Theta_{{\cal X}/S}^{\scriptboxtimes 2} &
  \rightarrow & \widetilde{\Theta}_{\cal X}
  \displayboxtimes\Theta_{{\cal X}/S} &
  \rightarrow
  \pi^{-1}\Theta_S\displayboxtimes\widetilde{\Theta}_{\cal X}
  \rightarrow(\pi^{-1}\Theta_S)^{\scriptboxtimes 2}
  \rightarrow 0 \\
 & \downarrow & & \downarrow & \\
 & \Theta_{{\cal X}/S} & = & \Theta_{{\cal X}/S} &
\end{array}
\end{equation}
where the vertical maps are composed of the restriction to
the diagonal $\Delta\subset{\cal X}\times_S{\cal X}$ followed
by Lie brackets.

\rk Here we use the fact that the restriction of
the Lie bracket
$$[\ ,\ ]:\Theta_{\cal X}^{\scriptboxtimes 2}
\longrightarrow\Theta_{\cal X}$$
to
$\widetilde{\Theta}_{\cal X}\boxtimes\Theta_{{\cal X}/S}$ takes
values in $\Theta_{{\cal X}/S}$ (see \cite{BS}, and
also \cite{EV}).

We note that the first column of the diagram (18) constitutes
the complex $\widetilde{\cal K}\udot$ computing
$\widetilde{\bf T}^{(2)}_{{\cal X}/S}$. Let ${\cal L}\udot$
denote the complex
$$\begin{array}{ccc}
  {\scriptstyle -1} & & {\scriptstyle 0} \\
  \widetilde{\Theta}_{\cal X}\displayboxtimes
  \Theta_{{\cal X}/S} & \stackrel{\ell}{\longrightarrow} &
  (\pi^{-1}\Theta_S\displayboxtimes
  \widetilde{\Theta}_{\cal X})\oplus\Theta_{{\cal X}/S}
\end{array}$$
with $\ell=(\beta\boxtimes\alpha,[\ ,\ ])$\ .
Then we can rewrite (18) as a short exact sequence
of {\em complexes}\ \ on ${\cal X}\times_S{\cal X}$:
\begin{equation}
 0\longrightarrow\widetilde{\cal K}\udot\longrightarrow
 {\cal L}\udot\longrightarrow
 (\pi^{-1}\Theta_S)^{\scriptboxtimes 2}\longrightarrow 0
\end{equation}
The associated direct-image sequence yields
\begin{equation}
 \longrightarrow{\bf R}^0f_*{\cal L}\udot
\longrightarrow\Theta_S^{\otimes 2}
 \stackrel{\widetilde{\kappa}_2}{\longrightarrow}
 \widetilde{\bf T}^{(2)}_{{\cal X}/S}\longrightarrow\ \ .
\end{equation}
Tracing out the definition of a connecting
morphism  (bearing in mind that if $\zeta\in\Theta_{\cal X}$ is
a local lifting of $\partial/\partial t_k$, and $\theta({\bf t})$ is
any element of $\Theta_{{\cal X}/S}$, then
$[\zeta,\theta({\bf t})]=
\frac{\partial\theta({\bf t})}{\partial t_k}$)
shows that it is indeed the same as $\widetilde{\kappa}_2$ given by
the explicit Definition in coordinates given above. \qed

\rk The explicit construction above shows how the data, up to
second order, of the (first) Kodaira-Spencer mapping
$$\kappa:\Theta_{S,0}/{\bf m}_{S,0}^2\Theta_{S,0}\longrightarrow
{\bf T}^1_{{\cal X}/S,0}/{\bf m}_{S,0}^2{\bf T}^1_{{\cal X}/S,0}$$
determines the second Kodaira-Spencer class
$$\widetilde{\kappa}_2:\Theta_S^{\otimes 2}|_0=
\Theta_{S,0}^{\otimes 2}/{\bf m}_{S,0}\Theta_{S,0}^{\otimes 2}
\longrightarrow
\widetilde{\bf T}^{(2)}_X$$
(see (12)). Conversely, if $\zeta,\xi\in\Theta_S|_0$, and
$\widetilde{\kappa}_2(\zeta\otimes\xi)$ is represented by a
cocycle $(\widehat{\zeta}\times\widehat{\xi},\theta)\in
\check{C}^1(\widetilde{\cal K}\udot|_0)$, then we can choose
coordinates $\bf t$ on $S$ so that $\zeta=\partial/\partial t_k$,
$\xi=\partial/\partial t_{\ell}$, and the Kodaira-Spencer mapping
of the deformation in question is represented
in
${\bf T}^1_{{\cal X}/S,0}\otimes\Omega^1_{S,0}/
{\bf m}_{S,0}^2{\bf
T}^1_{{\cal X}/S,0}\otimes\Omega^1_{S,0}$
by the form
$$\sum_{\mu=1}^s(\theta_{\mu}^{(0)}+\sum_{\nu=1}^s
\theta_{\mu}^{(\nu)}t_k)dt_{\ell}
$$
with
$\theta_k^{(0)}=\widehat{\zeta}$,
$\theta_{\ell}^{(0)}=\widehat{\xi}$, and
$\theta_{\ell}^{(k)}=\theta$.

\section{Main results}

There is a natural composition map
\begin{equation}
  \Theta_S^{\otimes 2}\hookrightarrow
  \Theta_S\oplus\Theta_S^{\otimes 2}\surj T_S^{(2)}\ .
\end{equation}
However, this map is not
${\cal O}_S$-linear. For example, $x\otimes y-y\otimes x$ is
mapped to $[x,y]$, whereas for any $f\in{\cal O}_S$
$$f.(x\otimes y-y\otimes x)=(f.x)\otimes y- y\otimes(f.x) $$
is sent to $f.[x,y]-y(f).x$ \ . Nevertheless, (21) induces an
${\cal O}_S$-linear map
$$\Theta_S^{\otimes 2}\surj T_S^{(2)}/\Theta_S\ \
(\simeq Sym^2\Theta_S)\ .$$
The latter fits in a commutative square of ${\cal O}_S$-linear
maps obtained from (8),
\begin{equation}
\begin{array}{ccc}
       \Theta_S^{\otimes 2} &
         \stackrel{\overline{\widetilde{\kappa}}_2}{\longrightarrow} &
         \widetilde{\bf T}^{(2)}_{{\cal X}/S}/
         \, im\, (\kappa_1) \\
       \downsurj & & \downsurj \\
        T_S^{(2)}/\Theta_S &
          \stackrel{\overline{\kappa}_2}{\longrightarrow} &
          {\bf T}^{(2)}_{{\cal X}/S}/\,im\,(\kappa_1)
\end{array}
\end{equation}

\begin{thm}
 The second differential of the Archimedean period map
 $d^2\Psi$
factors through the
 diagonal of (22),
$$\overline{\widetilde{\kappa}}_2:
\Theta_S^{\otimes 2}\longrightarrow
\widetilde{\bf T}_{{\cal X}/S}^{(2)}/\,im\,(\kappa_1)\ . $$
\end{thm}

\pf  Since the statement deals with ${\cal O}_S$-linear
maps, it is enough to prove it pointwise, for each $t\in S$. It
suffices to restrict to $0\in S$. We need to show that
$d^2\Psi(y)=0$ for any $y\in \Theta_S^{\otimes 2}|_0$ with
$\widetilde{\kappa}_2(y)\in im\,(\kappa_1)$. The condition on $y$
implies that
$\widetilde{\kappa}_2(y)\in \widetilde{\bf T}^{(2)}_X$
can be represented by a cocycle of the form
\begin{equation}
(0,\theta)\in \check{C}^2(\Theta_X^{\scriptboxtimes 2}[1])\oplus
\check{C}^1(\Theta_X)=\check{C}^1(\widetilde{\cal K}\udot|_0)\ ,
\end{equation}
where $\theta$ is a {\em cocycle} \ in
$\check{C}^1(\Theta_X)$ representing $\kappa_1(\eta)$ for some
$\eta\in \Theta_S|_0$.

At this point we ``recall" two theorems from \cite{K}.
\begin{thm}[(5.4) in \cite{K}]

If \
$\kappa\in {\bf T}^1_{{\cal X}/S,0}\otimes\Omega_{S,0}^1/
{\bf m}_{S,0}^2{\bf T}^1_{{\cal X}/S,0}\otimes\Omega_{S,0}^1$
is represented by the form
$$\sum_{\mu=1}^s(\theta_{\mu}^{(0)}+\sum_{\nu=1}^s
\theta_{\mu}^{(\nu)}t_k)dt_{\ell}
$$
with
$\theta_k^{(0)}=\widehat{\zeta}$,
$\theta_{\ell}^{(0)}=\widehat{\xi}$, and
$\theta_{\ell}^{(k)}=\theta$,
then the second differential of the Archimedean period map
$$d^2\Psi(\frac{\partial}{\partial t_k}\otimes
\frac{\partial}{\partial t_{\ell}}):
H_{ar}\longrightarrow F_{ar}^{-2}/
H_{ar}+span\,\{\nabla_{\eta}|_0({\cal H}_{ar})\ |\
\eta\in\Theta_S\}$$
is induced by the map
\begin{eqnarray*}
\lefteqn{H_{ar}=H_{ar}^m\longrightarrow} \\
 & & {\bf H}^m(\Omega\udot_X\otimes B_{ar}/
F_{ar}^0(\Omega\udot_X\otimes B_{ar})+
span\,\{ \check{\boldpounds}_{\eta}|_0{\cal
F}_{ar}^0\ |\  \eta\in\Theta_S\})
\end{eqnarray*}
given on the cochain level by
\begin{eqnarray}
\nonumber
\lefteqn{\omega_{i_1,\ldots,i_q}.T^p=\omega_Q.T^p\mapsto}\\
 & &
\widehat{\zeta}_{i_{-1}i_0}\contr\widehat{\xi}_{i_0i_1}\contr
 \omega_Q.T^p-\widehat{\xi}_{i_0i_1}\contr
 \check{\varpounds}_{\widehat{\zeta}_{i_0i_1}}
 \omega_Q.T^{p+1}+
 \theta_{i_0i_1}\contr\omega_Q.T^p\ .
\end{eqnarray}
\end{thm}

\

\begin{thm}[(5.7) in \cite{K}]

$d^2\Psi$ on $\Theta_S^{\otimes 2}|_0$ is determined by
$$\kappa\in {\bf T}^1_{{\cal X}/S,0}\otimes\Omega_{S,0}^1/
{\bf m}_{S,0}^2{\bf T}^1_{{\cal X}/S,0}\otimes\Omega_{S,0}^1\
.$$
\end{thm}

Reading the two theorems in light of the Remark at the end of
Section~2,  Theorem 4 shows that $d^2\Psi(y)$ is determined by
$\widetilde{\kappa}_2(y)=[(0,\theta)]\in \widetilde{\bf
T}^{(2)}_X$,
and Theorem 3 says that $d^2\Psi(y)$ is induced on the
cochain level by the contraction  with the {\em cocycle} \
$\theta$ representing $\kappa_1(\eta)$. This contraction is
equivalent to $\nabla_{\eta}|_0$ modulo $H_{ar}$  (in fact, it is
none other than $d\Psi(\eta)$), and so we have proved that
$d^2\Psi(y)=0$. \qed

\begin{thm}
The graded version of the second differential of the
Archimedean period map \ $\overline{d^2\Psi}$, as well as the
second differential of the usual period map \ $d^2\Phi$ and the
second fundamental form of the VHS, \ \II, \ all factor through
$$\kappa_1^2:\Theta_S^{\otimes 2}\longrightarrow
({\bf T}^1_{{\cal X}/S})^{\otimes 2}\ ,$$
and thus depend on \ $\kappa_1$ only.
\end{thm}

\pf Again it suffices to restrict to $0\in S$.
Suppose
$$\widetilde{\kappa}_2(y)=[(\widehat{\zeta}\times
\widehat{\xi},\theta)]\in \widetilde{\bf T}^{(2)}_X$$
for some $y=\zeta\otimes\xi\in\Theta_S^{\otimes 2}|_0$.
Examining formula (24), we observe that the term involving
$\theta$ lies in $F_{ar}^{-1}$. Therefore,
$\overline{d^2\Psi}(y)$ depends only on
$\kappa_1^2(y):=\kappa_1(\zeta)\otimes\kappa_1(\xi)$.
This proves the Theorem for $\overline{d^2\Psi}$.

The statements for $d^2\Phi$ and II follow from this by Lemma 2
(b) and Proposition 1, respectively. \qed

Finally, all the maps in question are symmetric, and so we may
pass from $\kappa_1^2$ to $Sym^2\kappa_1$ and from
$\overline{\widetilde{\kappa}}_2$ to $\overline{\kappa}_2$
(see (8)). Referring to the following symmetrized version of
(14),
\begin{equation}
\begin{array}{cccccccc}
  0\rightarrow & \Theta_S & \rightarrow &
    {\bf T}^{(2)}_S &
    \rightarrow & Sym^2\,\Theta_S & \rightarrow  & 0 \\
   & \kappa_1 \downarrow & & \downarrow \kappa_2 &
   & \downarrow Sym^2\kappa_1 & & \\
  0\rightarrow & {\bf T}_{{\cal X}/S}^1 & \rightarrow &
    {\bf T}_{{\cal X}/S}^{(2)} & \rightarrow &
    Sym^2{\bf T}_{{\cal X}/S}^1 &
    \stackrel{o}{\rightarrow} & {\bf T}_{{\cal X}/S}^2\ ,
\end{array}
\end{equation}
we conclude  with

\begin{thm} $d^2\Psi$ factors through
$$\overline{\kappa}_2:Sym^2\Theta_S\longrightarrow
{\bf T}_{{\cal X}/S}^{(2)}/\,im\,(\kappa_1),$$
whereas \ $\overline{d^2\Psi}$, $d^2\Phi$ and \ $\II$ \ factor
through
$$Sym^2\kappa_1:Sym^2\Theta_S\longrightarrow
Sym^2{\bf T}_{{\cal X}/S}^1\ .$$
\end{thm}

\rk When the deformation is versal, i.e.
$im\,(\kappa_1)$ is all of ${\bf T}_{{\cal X}/S}^1$, there is
no difference between $\overline{\kappa}_2$ and $Sym^2\kappa_1$.

\section{The higher-order cases}

The definition of the second differential of the period map in
Section 1 easily generalizes to higher-order cases (cf.
\cite{K}).

All three papers mentioned in the introduction define ``tangent
spaces to the moduli" ${\bf T}_{{\cal X}/S}^{(n)}$ of all
orders $n$. However, these definitions seem more complicated
than in the case $n=2$.

Still, we should have a diagram analogous to (25),
\begin{equation}
\begin{array}{cccccccc}
  0\rightarrow & {\bf T}^{(n-1)}_S & \rightarrow &
    {\bf T}^{(n)}_S &
    \rightarrow & Sym^n\,\Theta_S & \rightarrow  & 0 \\
   & \kappa_{n-1} \downarrow & & \downarrow \kappa_n &
   & \downarrow Sym^n\kappa_1 & & \\
  0\rightarrow & {\bf T}_{{\cal X}/S}^{(n-1)} & \rightarrow &
    {\bf T}_{{\cal X}/S}^{(n)} & \rightarrow &
    Sym^n{\bf T}_{{\cal X}/S}^1 &
    \stackrel{o_n}{\rightarrow} & {\bf T}_{{\cal X}/S}^2\ ,
\end{array}
\end{equation}
where $o_n$ is the $n^{th}$ obstruction map, and we expect that
the $n^{th}$ differential of the Archimedean period map
$d^n\Psi$ factors through the $n^{th}$ Kodaira-Spencer map
$\kappa_n$ modulo the image of $\kappa_{n-1}$, whereas the
$n^{th}$ differential of the usual period map $d^n\Phi$ and the
$n^{th}$  fundamental form of the VHS I$n$I factor through
$Sym^n\kappa_1$.

\section{Appendix: Archimedean cohomology}

In this section we summarize what we need about  Archimedean
cohomology. For more information on this subject we refer to
\cite{Den}.

\begin{dfntn}
$B_{ar}={\bf C}[T,T^{-1}], \ \ {\bf L}={\bf C}[T^{-1}]$. \
$B_{ar}$ is filtered by the ${\bf L}$-submodules
$F^p=T^{-p}.{\bf L}$.
\end{dfntn}

Thus, if $X$ is a compact K\"{a}hler manifold,
$H^m(X)\otimes_{\bf C}B_{ar}$ receives the filtration
$F\udot_{ar}$ obtained as the tensor product of the Hodge
filtration on $H^m(X,{\bf C})$ and the filtration $F\udot$ on
$B_{ar}$. $F\udot_{ar}$ is a decreasing filtration with
infinitely many levels, and
$$Gr_{F_{ar}}^k\cong\bigoplus_{p+q=m}H^{p,q}.T^{p-k}\ .$$

\begin{dfntn}
The {\em Archimedean cohomology}\ of $X$ is
$$H^m_{ar}(X):=F_{ar}^0(H^m(X,{\bf C})\otimes B_{ar})\ .$$
\end{dfntn}

Consider the complex of sheaves
$\Omega_X\udot\otimes_{\bf C}B_{ar}$ with the differential
$${\bf d}(\omega.T^k):=d\omega.T^{k+1}\ .$$
This complex is also filtered by the tensor product of the
stupid filtration on $\Omega\udot_X$ and $F\udot$ on $B_{ar}$,
and we have
\begin{eqnarray*}
\lefteqn{H^m_{ar}(X)=F^0{\bf H}^m(X,\Omega\udot_X\otimes
B_{ar}) \cong} \\
 & & {\bf H}^m(X,F^0(\Omega\udot_X\otimes B_{ar}))
\cong
 {\bf H}^m(X,\Omega\udot_X)\otimes{\bf L}
 \cong H^m(X,{\bf C})\otimes{\bf L}\ .
\end{eqnarray*}
Note that ${\bf H}^m(X,\Omega\udot_X\otimes B_{ar})$ is a
complex infinite-dimensional Hodge structure (of weight $m$),
and $(\Omega\udot_X\otimes B_{ar},{\bf d})$ is a Hodge complex.
Hence
$$Gr_{F_{ar}}^k{\bf H}^m(X,\Omega\udot_X\otimes B_{ar})\cong
{\bf H}^m(X,Gr_{F_{ar}}^k(\Omega\udot_X\otimes B_{ar}))\ .$$
We will write boldface $\check{\bf D}$ for the differential in
the \v{C}ech cochain complex computing
${\bf H}^m(X,\Omega\udot_X\otimes B_{ar})$, and
$$\check{\boldpounds}_v:=\check{\bf D}v\contr +
v\contr\check{\bf D}$$ for the corresponding Lie derivative
with respect to a vector field $v$ on $X$.

These constructions extend without any difficulty to the
relative situation. In particular, given a flat family
$\pi:{\cal X}\longrightarrow S$ of compact K\"{a}hler
manifolds, the bundle
$${\cal H}\otimes B_{ar}={\bf R}^m\pi_*(\Omega\udot_{{\cal
X}/S}\otimes B_{ar})$$
is filtered by ${\cal F}\udot_{ar}$, and the Gauss-Manin
connection $\nabla$ extends to
$$\nabla_{ar}:{\cal H}\otimes B_{ar}\longrightarrow
{\cal H}\otimes B_{ar}\otimes\Omega_S^1\ ,$$
with the usual Griffiths' transversality property
$$\nabla_{ar}({\cal F}^p_{ar})\subset{\cal F}^{p-1}_{ar}\otimes
\Omega_S^1\ .$$

Specifically, if $x$ is a section of $\cal H$, then
$$\nabla_{ar}(x.T^p)=\nabla x.T^p\ .$$

The real difference arises when one examines the definition of
$\nabla_{ar}$ on the cochain level, due to the fact that $\bf
d$ increases the exponent at $T$.

\

\

\

\noindent Department of Mathematics\\Columbia University\\
New York, NY 10027

\

\noindent {\em Current address}:

\

\noindent Department of Mathematics\\Yeshiva University\\
500 W 185 St.\\New York, NY 10033\hfill{\em yk@yu1.yu.edu}

\end{document}